\documentclass[a4paper,12pt]{amsart}
\usepackage{graphicx}
\usepackage{amssymb}
\usepackage{epic,eepic}
\allowdisplaybreaks[1]
\long\def\rem#1{}
\title{A conjectured formula for Fully Packed Loop configurations in a triangle}
\author[P.~Zinn-Justin]{Paul~Zinn-Justin}
\address{Paul Zinn-Justin, LPTMS (CNRS, UMR 8626), Univ Paris-Sud, 91405 Orsay Cedex, France; and LPTHE (CNRS, UMR 7589), Univ Pierre et Marie Curie-Paris6, 75252 Paris Cedex, France.}
\email{pzinn\,@\,lpthe.jussieu.fr}
\thanks{PZJ was supported by
EU Marie Curie Research Training Networks
``ENRAGE'' MRTN-CT-2004-005616, ``ENIGMA'' MRT-CT-2004-5652,
ESF program ``MISGAM''
and ANR program ``GIMP'' ANR-05-BLAN-0029-01.}
\thanks{The author wants to thank R.~Langer for suggesting the form of the
involution, J.~Thapper for sharing his numerical data,
as well as P.~Di Francesco, T.~Fonseca and J.-B.~Zuber for discussions.}
\numberwithin{equation}{section}
\newtheorem{prop}{Proposition} 
\newtheorem{lemma}{Lemma}
\newtheorem{conj}{Conjecture}
\newtheorem*{conj*}{Conjecture}
\date{June 2008}

\newcommand\A{{\mathbf{A}_n}}
\newcommand\Ap{{\mathbf{A}_{n+1}}}
\newcommand\Am{{\mathbf{A}_{n+m}}}
\newcommand\alg{{\Lambda_n}}
\newcommand\alginf{{\Lambda_\infty}}
\newcommand\alginfhat{{\hat{\Lambda}_\infty}}
\newcommand\zeron{\varnothing}
\newcommand\onea{\begin{picture}(6,8)
\put(0,0){\line(0,1){8}}
\put(0,8){\line(1,0){8}}
\put(0,0){\line(1,1){8}}
\end{picture}}
\newcommand\oneb{\begin{picture}(4,5)
\put(0,0){\line(0,1){5}}
\put(0,5){\line(1,0){5}}
\put(0,0){\line(1,1){5}}
\end{picture}}
\newcommand\one{\mathchoice{\onea}{\onea}{\oneb}{\oneb}}
\newcommand\onen{{\one_n}}
\newcommand\onenp{{\one_{n+1}}}
\newcommand\inv{\imath}
\newcommand\matinv{\mathcal{I}}
\newcommand\LR{\mathcal{C}}
\newcommand\LRt{\mathcal{\tilde C}}
\newcommand\braket[1]{\left\langle#1\right\rangle}
\renewcommand\t{\boldsymbol{\tau}}
\newdimen{\cellsize}

\newcommand\smallboxes{\setlength{\cellsize}{7pt}\def\boxformat{\scriptstyle}}
\def\boxformat{}
\newsavebox{\cellcontent}
\def\hidehrule#1#2{\kern-#1
  \hrule height#1 depth#2 \kern-#2 }%
\def\hidevrule#1#2{\kern-#1{\dimen\cellcontent=#1%
    \advance\dimen\cellcontent by#2\vrule width\dimen\cellcontent}\kern-#2 }%
\def\makeblankbox#1#2{\hbox{\lower\dp\cellcontent\vbox{\hidehrule{#1}{#2}%
    \kern-#1 
    \hbox to \wd\cellcontent{\hidevrule{#1}{#2}%
      \raise\ht\cellcontent\vbox to #1{}
      \lower\dp\cellcontent\vtop to #1{}
      \hfil\hidevrule{#2}{#1}}%
    \kern-#1\hidehrule{#2}{#1}}}}
\newcommand\cellify[1]{\defaultcell%
\sbox{\cellcontent}{\vbox to \cellsize{%
\vfill%
\hbox to \cellsize{\hfill$\boxformat #1$\hfill}
\vfill}}%
\rlap{\drawnbox}
\usebox{\cellcontent}}
\newcommand\tableau[1]{\vtop{\let\\\cr
\baselineskip -16000pt \lineskiplimit 16000pt \lineskip 0pt
\ialign{&\cellify{##}\cr#1\crcr}}}
\newcommand\defaultcell{\gdef\drawnbox{
\makeblankbox{0.2pt}{0.2pt}
}}
\newcommand\vdotscell{\gdef\drawnbox{\kern-1.6pt\vbox{\baselineskip=4pt\lineskiplimit=0pt\hbox{}\hbox{.}\hbox{.}\hbox{.}\hbox{}}}}
\newcommand\hdotscell{\gdef\drawnbox{\vbox to \cellsize{\hbox{\kern1pt$\ldotp\ldotp\ldotp$}}}}
\newcommand\vhdotscell{\gdef\drawnbox{\rlap{\kern-1.6pt\vbox{\baselineskip=4pt\lineskiplimit=0pt\hbox{}\hbox{.}\hbox{.}\hbox{.}\hbox{}}}\vbox to \cellsize{\hbox{\kern1pt$\ldotp\ldotp\ldotp$}}}}
\begin{document}
\begin{abstract}
We describe a new conjecture involving Fully Packed Loop counting
which relates recent observations of Thapper to formulae in the
Temperley--Lieb model of loops, and how it implies the Razumov--Stroganov
conjecture.
\end{abstract}
\maketitle
\section{Introduction}
In the literature generated by the seminal papers 
\cite{BdGN-XXZ-ASM-PP,RS-conj} 
and revolving around the so-called Razumov--Stroganov (RS) conjecture, 
it has often been remarked
that there are more conjectures than theorems. The present work, sadly,
will not help correct this imbalance: it is entirely based on one more
conjecture. The latter, however, is of some interest since it connects
the two sides of the Razumov--Stroganov conjecture; that is, it
expresses the number of Fully Packed Loop configurations
(FPLs) in a triangle with certain boundary conditions
as the constant term of a quasi-generating function which is closely
related to expressions appearing 
in the context of the Temperley--Lieb$(1)$ (sometimes
called $O(1)$) model of loops \cite{artic41,artic42,artic43}. 
This formula was inspired by an attempt
to understand the observations of Thapper \cite{Thapper} on the enumeration
of FPLs with prescribed connectivity, itself based on earlier
work \cite{CK,CKLN}. In fact, we shall show in what follows that our new
conjecture implies both the RS conjecture \cite{RS-conj} 
and the conjectures of \cite{Thapper}.

The paper is organized as follows. 
The next section contains some basic definitions.
In section 3, we briefly recall the various statistical models involved and the required content
from the work referred to above.
Section 4 contains the main formula of this paper, its conjectural meaning,
and the connection to the Razumov--Stroganov conjecture. 
Section 5 provides the link to Thapper's conjectures. 
The final section briefly describes the introduction of an extra parameter
$\t$ in the formulae.
Technical details
and numerical results are to be found in the appendices.

\section{Preliminaries}
\subsection{Bijections}\label{secbij}
\begin{figure}[t]
\caption{Bijections. From left to right: Ferrers diagrams and sequences of increasing integers; Ferrers diagrams and Dyck paths; Dyck paths and link patterns.}\label{figbij}
\includegraphics[scale=0.65]{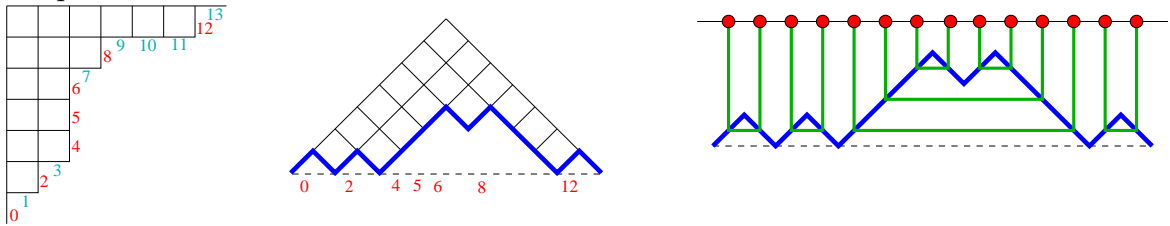}
\end{figure}

Let $n$ be a positive integer.
Various sets are in bijection:
\begin{enumerate}
\item
the set of Ferrers diagrams 
contained inside the ``staircase'' 
Ferrers diagram
(denoted by $\onen$ in what follows) with rows of lengths $n-1,n-2,\ldots,1$; 
\item
the set of Dyck paths of length $2n$;
\item
the set of link patterns of size $2n$, that is planar pairings of $2n$ points inside a half-plane (the points
sitting on the boundary); 
\item
the set of sequences of integers $(\alpha_i)_{i=0,\ldots,n-1}$
such that $\alpha_{i+1}>\alpha_i$ and $0\le \alpha_i\le 2i$ for all $i$.
\end{enumerate}
These various bijections are described on Fig.~\ref{figbij}. The only bijection
which we shall write down explicitly is from Ferrers diagrams to increasing
sequences: starting from a diagram $\alpha\subset\onen$,
consider the sequence of lengths of its rows and
pad it with zeroes so that it has exactly $n$ parts 
$\tilde\alpha_1,\ldots,\tilde\alpha_n$
(with in particular $\tilde\alpha_n=0$); the sequence is then given by
$\alpha_i=\tilde\alpha_{n-i}+i$, $i=0,\ldots,n-1$. 

We shall mostly use in what follows Ferrers diagrams and increasing
sequences, identified via the bijection above.
We shall call $\A$ their set.

\subsection{Order, embedding}\label{secembed}
We call $|\alpha|$ the number of boxes of 
$\alpha\in \A$: $|\alpha|=\sum_{i=0}^{n-1} (\alpha_i-i)$.

We consider the partial order on Ferrers diagrams which is simply
inclusion. In terms of sequences, $\alpha\subset \beta$ iff $\alpha_i\le\beta_i$
for all $i$. The smallest element is the empty Ferrers diagram, denoted by
$\zeron$; the largest element is $\onen$ itself.

In what follows whenever we consider matrices with indices
in $\A$, we shall assume that an arbitrary total order which
refines $\subset$ has been chosen. Accordingly,
an upper triangular matrix $M_{\alpha\beta}$ is
a matrix such that $M_{\alpha\beta}=0$ whenever $\alpha\not\subset \beta$.

Finally, note that viewed as Ferrers diagrams,
we have $\A\subset\Ap$.
This embedding, in terms of sequences,
sends $\alpha=(\alpha_0,\ldots,\alpha_{n-1})$ to
$(0,\alpha_0+1,\ldots,\alpha_{n-1}+1)$.
One must be warned that not all quantities defined below
satisfy a ``stability'' property with respect to this embedding
i.e.\ some quantities depend explicitly on $n$ and not just on the
underlying Ferrers diagrams. The quantities that are stable are
the matrices $\mathbf{P}$, $\mathbf{C}$, $\LR$, $\matinv$ to be defined below.
The quantities that are not stable satisfy instead
recurrence relations with respect to $n$, see section~\ref{secrecur}.

\subsection{Schur functions}
To any Ferrers diagrams one can associate a Schur function.
In the case of $\alpha\in\A$, and using an alphabet of $n$ letters,
the Schur function can be defined in terms of 
the corresponding sequence of integers $\alpha=(\alpha_i)$ as
\begin{equation}\label{defschur}
s_\alpha(u)=\frac{\displaystyle\det\left(u_i^{\alpha_j}\right)_{0\le i,j\le n-1}}
{\Delta(u)}
\end{equation}
where $u:=(u_0,\ldots,u_{n-1})$ and the denominator
is simply the numerator with empty Ferrers diagram:
$\Delta(u):=\prod_{0\le i<j\le n-1}(u_j-u_i)$, 
that is the Vandermonde determinant.


\subsection{The involution}
The Schur functions associated to $\alpha\in\A$ span
a vector space of dimension $c_n$. It can be made into
a commutative algebra $\alg$ by defining its structure constants to be
the Littlewood--Richardson coefficients $\LR_{\sigma\tau}^{\rho}$.
In other words $\alg$ is a truncation of the algebra of symmetric
functions in which one restricts oneself to diagrams inside
$\onen=(n-1,n-2,\ldots,1)$: if $\alginf$ denotes the algebra of
symmetric functions, that is simply the algebra of all Schur functions
with the ordinary function product, then
$\alg$ is canonically identified with the quotient
$\alginf$ by the span of the $\sigma\not\subset \onen$, the latter
being an ideal. In what follows we actually need the slightly larger
space $\alginfhat$ of symmetric power series. $\alg$ is clearly
also a quotient of $\alginfhat$.

We now introduce an involution $\inv$ on $\alginfhat$.
It is defined by its action on elementary symmetric functions $e_i$
through their generating series $\prod_i (1+z u_i)=\sum_i e_i z^i$:
\begin{equation}\label{defphi}
\inv\left(\prod_i (1+z u_i)\right)=\prod_i \frac{1}{1-z\frac{u_i}{1+u_i}}
\end{equation}
where the $(u_i)$ are an arbitrary alphabet,
and the equality should be understood order by order in $z$.
By the morphism property this defines $\inv$ entirely on $\alginf$
(we shall extend it below to $\alginfhat$).

Denote
\begin{equation*}
\tilde s_\alpha(u):=\inv(s_\alpha(u))
\end{equation*}
One can compute $\tilde s_\alpha(u)$ explicitly as follows.
Note that by \eqref{defphi}, $\inv$ is the composition of:
(a) the change of variables $u\mapsto
\frac{u}{1+u}$, and (b) the transposition of diagrams.\footnote{Equivalently,
it is the composition of two commuting involutions: 
(a) $u\mapsto -\frac{u}{1+u}$ and (b) transposition
of Ferrers diagrams
composed with multiplication by $(-1)^{|\alpha|}$.}
Thus, if one defines the sequence
$(\alpha'_i)$ associated to the transpose diagram
$\alpha'$ as the ordered complement of the $\{ 2n-1-\alpha_i \}$
inside $\{0,1,\ldots,2n-1\}$
(it can also be defined from the lengths of the columns $\tilde\alpha'_i$ by
$\alpha'_i=\tilde\alpha'_{n-i}+i$), then from \eqref{defschur},
\begin{equation}\label{tildas}
\begin{split}
\tilde s_\alpha(u)&=\frac{\det\left(\left(\frac{u_i}{1+u_i}\right)^{\alpha'_j}\right)_{0\le i,j\le n-1}}{\Delta\left(\frac{u}{1+u}\right)}
\\&=\prod_{i=0}^{n-1} (1+u_i)^{n-1}
\frac{\det\left(\left(\frac{u_i}{1+u_i}\right)^{\alpha'_j}\right)_{0\le i,j\le n-1}}{\Delta(u)}
\end{split}
\end{equation}

This also leads to the following lemma:
\begin{lemma}\label{exptilde}
There is an expansion of the form
\begin{equation*}
\tilde s_\alpha(u)=s_{\alpha'}(u)+\sum_{\beta\supsetneq\alpha'} c_\beta s_{\beta}(u)
\end{equation*}
where the $c_\beta$ are some coefficients.
\end{lemma}
\begin{proof}
Expand by multilinearity 
\[
\frac{\det\left(\left(\frac{u_i}{1+u_i}\right)^{\alpha'_j}\right)_{0\le i,j\le n-1}}{\Delta(u)}
=
\sum_{k_1,\ldots,k_n\ge 0} 
c_{k_1,\ldots,k_n}
\frac{\det\left(u_i^{\alpha'_j+k_j}\right)_{0\le i,j\le n-1}}{\Delta(u)}
\]
where the $c_{k_1,\ldots,k_n}$ are irrelevant binomial coefficients,
with $c_{0,\ldots,0}=1$.

Note that the sequence $(\alpha'_i+k_i)_{i=0,\ldots,n-1}$
is not necessarily increase;
however, if two terms are equal, then the determinant is zero,
and if they are all distinct, then there exists a permutation $\mathcal{P}$
which reorders them; call $\beta$ the corresponding increasing
sequence: $\beta_{\mathcal{P}(i)}:=\alpha'_i+k_i$.
In the latter case we have $\det(u_i^{\beta_j})/\Delta(u)
=(-1)^{|\mathcal{P}|} s_\beta(u)$.

Next we use the fact, which will be needed again in what follows, that if 
$\alpha'_i\le\beta_{\mathcal{P}(i)}$ for some $\mathcal{P}$ and all $i$,
where $\alpha'$ and $\beta$
are increasing sequences, then $\alpha'_i\le\beta_i$ for all $i$
(induction on the number of inversions
$|\mathcal{P}|$, noting that if $i$ is such that
$\mathcal{P}(i)>\mathcal{P}(i+1)$, then
$\alpha'_i<\alpha'_{i+1}\le \beta_{\mathcal{P}(i+1)}<\beta_{\mathcal{P}(i)}$, so that
one can permute the images of $i$ and $i+1$, thus reducing $|\mathcal{P}|$ by one).

Combining the facts above leads to the expansion of the form of the lemma.
\end{proof}

This lemma has two consequences. The first is that 
$\inv$ is well-defined on the whole of $\alginfhat$
(only finite sums occur for any coefficient of the image
of any symmetric power series).
It is then easy to
check that $\inv$ is indeed an involution on $\alginfhat$.

The second is that this involution $\inv$ is compatible with
the quotient to $\alg$
(keeping in mind that $\onen$ is invariant by transposition).

Finally, note that setting $z=1$ in \eqref{defphi} leads to
\begin{equation}\label{selfdual}
\inv\left(\prod_i (1+ u_i)\right)=\prod_i (1+u_i)
\end{equation}
Thus, $\prod_i (1+u_i)$, the sum of elementary symmetric functions,
is left invariant by $\inv$.

\subsection{Change of basis}
The link patterns can be considered as forming the canonical basis of a vector space.
It is however convenient to introduce another
basis; it is 
naturally indexed by elements of $\A$ (increasing sequences) too, 
and is related to the canonical basis by a triangular change of basis
(recall that we identify indices using the bijection of \ref{secbij}).
If a vector has entries $\psi_\pi$ in the canonical basis and entries
$\Psi_\alpha$ in this new basis
(note the use of lowercase vs uppercase in order to distinguish these quantities), then
\begin{equation}\label{chgbasis}
\Psi_\alpha=\sum_{\pi\in \A} 
\psi_{\pi} \mathbf{P}^\pi_{\ \alpha} 
\end{equation}
Here $\mathbf{P}$ is the transpose of the usual matrix of change of
basis. The reason for this transposition 
is that, to conform with the conventions of \cite{Thapper},
our operators will act on the right.
The matrix $\mathbf{P}$ is given explicitly in appendix \ref{Pmat};
we only need the fact that it is upper triangular, with ones on the diagonal.

\section{Statistical models of loops and Razumov--Stroganov conjecture}
\subsection{Counting of FPLs}\label{secfpl}
We introduce here the statistical lattice model called Fully Packed Loop (FPL) model.
It is defined on a subset of the square lattice;
in any given configuration of the model, edges of the lattice can have two states, empty or occupied,
in such a way that each vertex has exactly two neighboring occupied edges (i.e.\ paths made of occupied edges visit every vertex of the
lattice). We only consider the situation in which the Boltzmann weights are trivial, that is the pure enumeration problem.

\begin{figure}
\includegraphics[width=5cm]{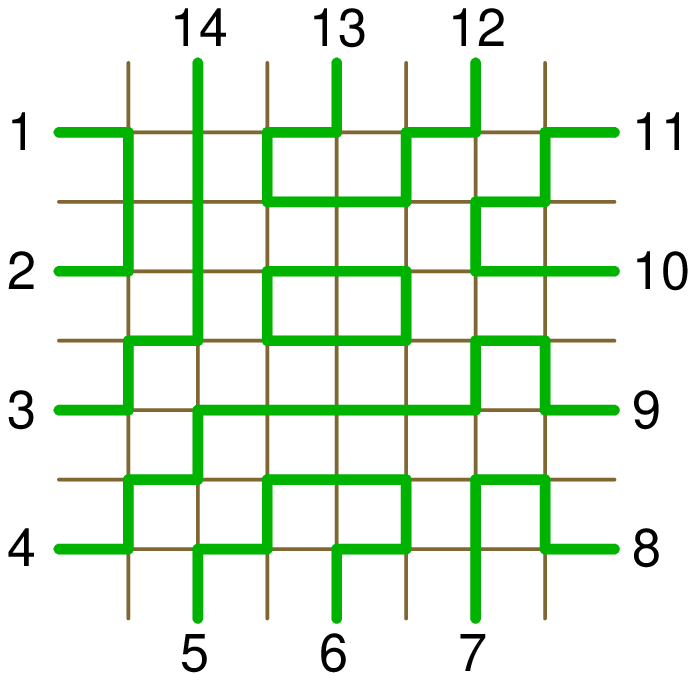}\qquad
\includegraphics{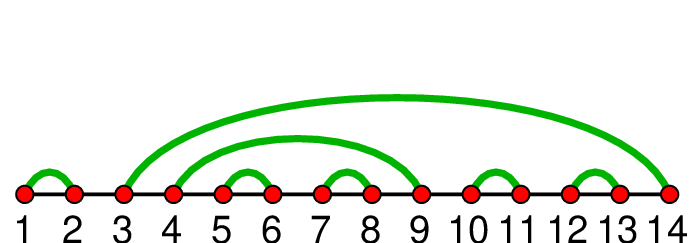}
\caption{A FPL configuration and the associated link pattern.}\label{fpl-ex}
\end{figure}

Given a positive integer $n$, we are interested in FPL configurations inside a $n\times n$ grid with specific boundary conditions
exemplified on Fig.~\ref{fpl-ex}: external edges are alternatingly occupied or empty. The justification for these boundary conditions
comes from the connection to the six-vertex model (in which they correspond to the so-called Domain Wall Boundary Conditions),
as well as to Alternating Sign Matrices, see \cite{Propp-ASM}.

Observe that there are two types of paths: the closed ones (loops) and the open ones, whose endpoints lie on the boundary.
Ignoring the former, we see that to each FPL we can associate a link pattern that encodes the connectivity of its 
endpoints. Let us call $\psi_\pi$ the number of FPLs with connectivity $\pi$.
Note that the endpoints must be labelled, which implies the choice of an origin;
but it is in fact irrelevant due to Wieland's theorem \cite{Wieland}, which states that $\psi_\pi = \psi_{\rho(\pi)}$
where $\rho(\pi)$ is the link pattern obtained from $\pi$ by cyclic rotation of the $2n$ points.

\begin{figure}
\includegraphics[width=7cm]{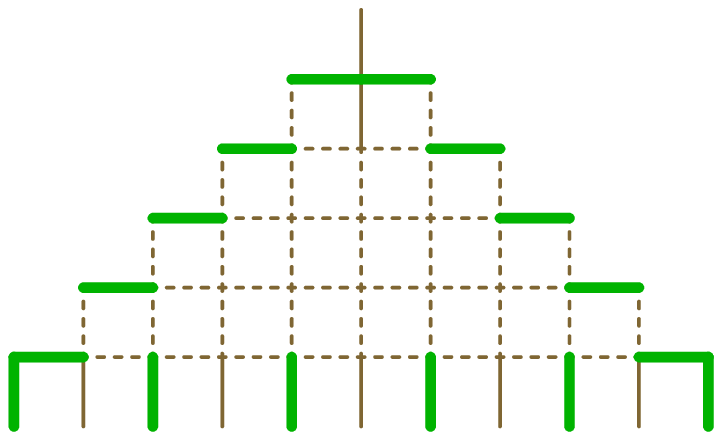}
\caption{Boundary conditions for FPLs in a triangle.}\label{trifpl-template}
\end{figure}

\smallboxes
\newcommand\vcenterbox[1]{\vcenter{\hbox{#1}}}
\begin{figure}
$\vcenterbox{\includegraphics[width=7cm]{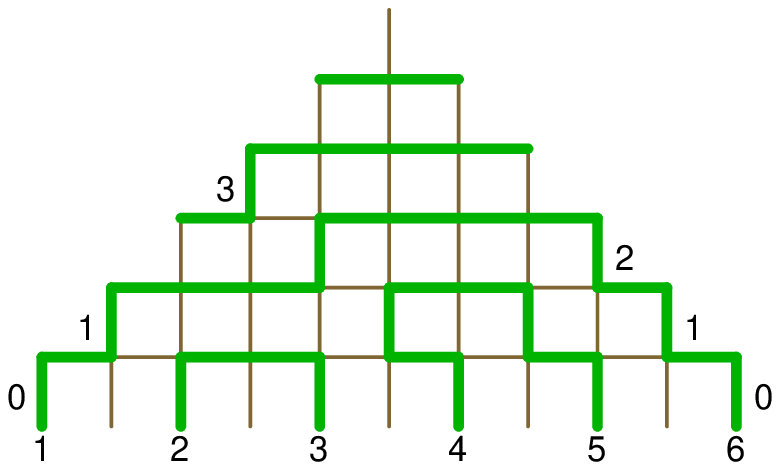}}\qquad
\begin{aligned}
\pi&=\includegraphics[width=4.5cm]{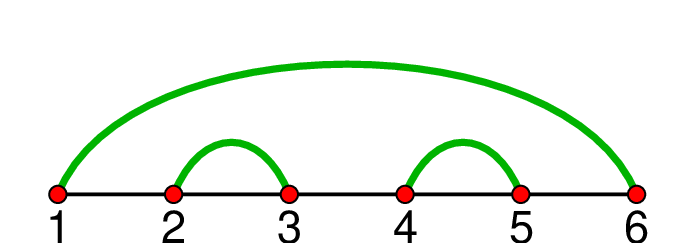}\\
\sigma&=(0,1,3)=\tableau{\\}\\
\tau&=(0,1,2)=\zeron
\end{aligned}
$
\caption{Example of parameterization of the boundary conditions of FPLs in a triangle.}\label{trifpl-ex}
\end{figure}

In general, one does not know how to compute $\psi_\pi$. There has been however
some progress \cite{dG-review,artic27,artic30,CKLN,Thapper}, and we are particularly
interested here in the formulae of \cite{CKLN,Thapper}, which appear
as a byproduct of proofs or attempted proofs of certain conjectures of \cite{Zuber-conj}.
Specifically, consider as in \cite{Thapper} FPL configurations in a triangle
of the form of Fig.~\ref{trifpl-template}, with exactly $2n$ vertical occupied external edges at the bottom, interlaced
with $2n-1$ empty edges. We further
require that each of the $2n-2$ external horizontal edges on the left (excluding the bottommost one) 
be connected to one of the $2n-2$ external horizontal edges on the right, and vice versa.
These configurations can
be classified as follows: the connectivity of the vertical external edges
can be encoded into a link pattern $\pi$ of size $2n$; furthermore, it is shown
in \cite{Thapper} 
that if one considers the sequence of the $2n$ vertical edges on either left or right boundaries
read from bottom to top, then it forms a Dyck path with occupied=up and empty=down. 
Equivalently, in our language, the sequence of locations of occupied vertical edges on the
left (resp.\ right) boundary,
numbered from bottom ($0$) to top ($2n-1$), is a sequence in $\A$, which we 
denote by $\sigma$ (resp.\ $\tau$), see Fig.~\ref{trifpl-ex} for an example.
Finally, define
$a_{\sigma,\pi,\tau}$ to be the number of FPLs in a triangle with the boundary conditions
defined above and given $\sigma$, $\pi$, $\tau$.

\rem{important to understand what is really constrained or not in these configurations
if one hopes to prove my conjecture}

Then the following equality between the two enumeration problems holds:
\begin{equation}\label{presummation}
\psi_\pi:=\sum_{\sigma,\tau\in \A} a_{\sigma,\pi,\tau} P_{\sigma'}(-k) P_{\tau'}(k-n+1)
\end{equation}
where $k$ is an integer to be discussed below, and
$P_\sigma(x)$ is a polynomial of $x$ which coincides for $x$ positive
integer with the number of SSYT of shape $\sigma$; in fact,
\begin{equation*}
P_\sigma(x)=
\begin{cases}
s_\sigma(\underbrace{1,\ldots,1}_x)&x\in \mathbb{Z}_+\\
s_{\sigma'}(\underbrace{-1,\ldots,-1}_{-x})
&x\in\mathbb{Z}_-
\end{cases}
\end{equation*}
Explicitly, it is given by
\begin{equation*}
P_\sigma(x)=
\prod_{(i,j)\in\lambda}\frac{j-i+x}{\lambda_i-i+\lambda'_j-j+1}
\end{equation*}
where $\lambda_i$ (resp.\ $\lambda'_i$) is the length of the $i^{\mathrm{th}}$
row (resp.\ column) of $\sigma$.

Formula \eqref{presummation} can be
deduced from Eq.~(4) in \cite{Thapper}.\footnote{Technically,
it is obtained from Eq.~(4) of \cite{Thapper} by setting $m=0$ in it,
$m$ being the number of extra arches that surround all arches (see also
section \ref{secrecur} of the present work).
The formula of Theorem 4.2 is only proved for $m$ sufficiently large,
but a clever argument in section 5 of \cite{CKLN} shows a property of polynomiality in $m$,
which allows to continue it to $m=0$.}
In the derivation, the value of $k$ appears in relation to the geometry
of FPLs, and the exact range of $k$ for which the formula is proved is not made clear.
In the present work we only require that the formula be true for one value of $k$ -- 
the explicit value being irrelevant since
the result should be independent of $k$. See also Theorem 4.2 of \cite{CKLN} (which is the
special case $k=0$).

\subsection{Temperley--Lieb$(1)$ loop model}
Another, {\em a priori}\, unrelated model is the following.
Consider the semi-group generators $\mathbf{e}_i$, $1\le i\le 2n$, acting on link patterns of size $2n$
as follows: $\mathbf{e}_i$ turns a link pattern $\pi$ into the link pattern obtained from $\pi$
by pairing together (i) the points which are connected to $i$ and $i+1$ in $\pi$ and (ii) $i$ and $i+1$,
all the other pairings remaining the same. For $i=2n$ one assumes periodic boundary conditions i.e.\ $2n+1\equiv 1$.
By linearity the $\mathbf{e}_i$ can be made into operators on the space of linear combinations of link patterns
(thus forming a representation of the Temperley--Lieb$(1)$ algebra, see \cite{dG-review} for more details)
and one can then define the {\em Hamiltonian}:
\[
\mathbf{H}=\sum_{i=1}^{2n} \mathbf{e}_i
\]
The $\mathbf{e}_i$ (resp.\ $\mathbf{H}$) 
possess a left eigenvector which is $(1,\ldots,1)$ in the canonical
basis, with eigenvalue $1$ (resp. $2n$); thus, $\mathbf{H}$ also possesses a (right) eigenvector with the
same eigenvalue, denoted by $\psi'$:
\begin{equation}\label{eigenvec}
\mathbf{H} \psi' = 2n \psi'
\end{equation}
It is easy to check that $\mathbf{H}$ satisfies the hypotheses of the Perron--Frobenius theorem, 
so that $2n$ is the (strict) largest eigenvalue of $\mathbf{H}$, and $\psi'$ is uniquely defined by \eqref{eigenvec}
up to normalization. The latter,
since $\mathbf{H}$ has integer entries, can always be chosen such that $\psi'$ has positive coprime integer entries,
denoted by $\psi'_\pi$. An example is provided in appendix \ref{PFTL}.

In a series of papers \cite{artic31,artic34,artic41}, it was shown how to generalize the
Temperley--Lieb$(1)$ loop model to an inhomogeneous model, then relate its Perron--Frobenius eigenvector to the
quantum Knizhnik--Zamolodchikov equation, and finally write quasi-generating functions for entries of $\psi'$.
More precisely, the last step involves first performing the change of basis \eqref{chgbasis}; then
the new entries $\Psi'_\alpha$ can be written
\begin{equation}\label{prepsiformula}
\Psi'_\alpha=
\Delta(u)
\prod_{0\le i<j\le n-1}(1+u_j+u_i u_j)\ \Big|_{\prod_{i=0}^{n-1} u_i^{\alpha_i}}
\end{equation}
where $\Delta(u)=\prod_{0\le i<j\le n-1}(u_j-u_i)$, and $|_{\cdots}$ means picking the coefficient of a monomial
in a polynomial of the variables $u_0,\ldots,u_{n-1}$.
(note that in \cite{artic41} a slightly different notation $a_i=1+\alpha_{i-1}$, $1\le i\le n$ is used).

\subsection{Razumov--Stroganov conjecture}
Finally, to conclude this introductive section, we mention the following conjecture as main motivation for this work:
\begin{conj*} (Razumov, Stroganov \cite{RS-conj})
For $\pi$ a link pattern of size $2n$, let $\psi'_\pi$ be as above the entry of the properly normalized 
Perron--Frobenius eigenvector of the Hamiltonian of the Temperley--Lieb(1) loop model,
and $\psi_\pi$ be the number of FPLs with connectivity $\pi$; then
\[
\psi'_\pi = \psi_\pi 
\]
\end{conj*}

\section{The formula}
Let $n$ be a fixed positive integer,
$\sigma$, $\tau$ be two Ferrers diagrams
and $\alpha=(\alpha_i)$ be a sequence 
of integers in $\A$. We define $A_{\sigma,\alpha,\tau}$ to be the
coefficient of a monomial in the expansion of a certain formal power series:
\begin{equation}\label{mainformula}
A_{\sigma,\alpha,\tau}=
\tilde s_\sigma(u) s_\tau(u) 
\Delta(u)
\prod_{i=0}^{n-1} (1+u_i)^{n-1}
\prod_{0\le i<j\le n-1}(1+u_j+u_i u_j)\ \Big|_{\prod_{i=0}^{n-1} u_i^{\alpha_i}}
\end{equation}

Note the important fact 
that $\prod_{0\le i<j\le n-1}(1+u_j+u_i u_j)$ is a nonsymmetric factor.
If it were symmetric, we would
simply be picking one term in the expansion 
of a certain symmetric function in terms of Schur functions, but it is not so.

One can rewrite \eqref{mainformula} as a multiple contour integral in which the
contours surround $0$ clockwise (but not $-1$):
\begin{equation}\label{intformula}
A_{\sigma,\alpha,\tau}=\oint \prod_{i=0}^{n-1} \frac{d u_i}{2\pi i u_i^{\alpha_i+1}}
\tilde s_\sigma(u) s_\tau(u) 
\Delta(u)
\prod_{i=0}^{n-1} (1+u_i)^{n-1}
\prod_{0\le i<j\le n-1}(1+u_j+u_i u_j)
\end{equation}

Using \eqref{defschur} and \eqref{tildas}, one can also rewrite \eqref{mainformula}
more explicitly:
\begin{multline}\label{intformulab}
A_{\sigma,\alpha,\tau}=
\frac{\det\left(u_i^{\tau_j}\right)
\det\left(\left(\frac{u_i}{1+u_i}\right)^{\sigma'_j}\right)
}{\Delta(u)} 
\prod_{i=0}^{n-1} (1+u_i)^{2(n-1)}
\prod_{0\le i<j\le n-1}(1+u_j+u_i u_j)
\ \Big|_{\prod_{i=0}^{n-1} u_i^{\alpha_i}}
\end{multline}
where the $\tau_j$ and $\sigma'_j$ are the increasing sequences
associated to $\tau$ and $\sigma'$.

In what follows we shall use the simplifying notation: let us write
\begin{equation}
\braket{F(u)}_\alpha:=
F(u)\Delta(u)
\prod_{0\le i<j\le n-1}(1+u_j+u_i u_j)\ \Big|_{\prod_{i=0}^{n-1} u_i^{\alpha_i}}
\end{equation}
for any symmetric function $F$. With this notation,
\begin{equation}\label{defA}
A_{\sigma,\alpha,\tau}=\braket{\tilde s_\sigma(u) s_\tau(u) 
\prod_{i=0}^{n-1} (1+u_i)^{n-1}}_\alpha
\end{equation}
$A_{\sigma,\alpha,\tau}$ is of course an integer.

\subsection{Some properties}
An important fact is the following:
\begin{lemma}\label{nadeauthmb}
If $|\sigma|+|\tau|>|\alpha|$, $A_{\sigma,\alpha,\tau}=0$.
If $|\sigma|+|\tau|=|\alpha|$, $A_{\sigma,\alpha,\tau}=\LR_{\sigma',\tau}^{\alpha}$.
\end{lemma}
\begin{proof}
By degree counting. According to
lemma \ref{exptilde}, the lowest degree terms of \eqref{mainformula}
as a power series in the variables $u_i$ are
$s_{\sigma'}(u) s_\tau(u) \Delta(u)$.
i.e.\ of degree $|\sigma|+|\tau|+n(n-1)/2$. The first result follows from the fact
that we pick out a term of degree $\sum_i \alpha_i=|\alpha|+n(n-1)/2$.
If $|\sigma|+|\tau|=|\alpha|$,  one finds
$A_{\sigma,\alpha,\tau}=s_{\sigma'}(u) s_\tau(u)\Delta(u)_{|\prod_i u_i^{\alpha_i}}$.
Expanding $s_{\sigma'}(u) s_\tau(u)\Delta(u)=
\sum_\rho\LR_{\sigma',\tau}^\rho s_\rho(u)\Delta(u)=
\sum_\rho\LR_{\sigma',\tau}^\rho \det(u_i^{\rho_j})$, we conclude that
$A_{\sigma,\alpha,\tau}=\LR_{\sigma',\tau}^{\alpha}$.
\rem{equivalently one could use the scalar product on schur functions
as in the proof right below}
\end{proof}

Compare the first part of the lemma with Lemma 3.7 of \cite{Thapper}.
By triangularity of the change of basis $\mathbf{P}$,
the second part of the lemma also says that
$a_{\sigma,\alpha,\tau}=\LR_{\sigma',\tau}^{\alpha}$.
This generalizes lemma 3.6(b) of \cite{Thapper}, that
is lemma 4.1(2) of \cite{CKLN}.

Similarly, we have
\begin{lemma} \label{uptri}
If $\tau\not\subset \alpha$ or $\sigma'\not\subset\alpha$, 
$A_{\sigma,\alpha,\tau}=0$.
\end{lemma}
\begin{proof}
By symmetrizing the integrand of \eqref{intformula}, one can write
\[
A_{\sigma,\alpha,\tau}=\oint \prod_{i=0}^{n-1} \frac{d u_i}{2\pi i u_i}
\tilde s_\sigma(u) s_\tau(u) 
\Delta(u)
\prod_{i=0}^{n-1} (1+u_i)^{n-1}
\mathrm{AS}\left[
\frac{\prod_{0\le i<j\le n-1}(1+u_j+u_i u_j)}{\prod_{i=0}^{n-1} u_i^{\alpha_i}}
\right]_{\le0}
\]
where $\mathrm{AS}(f):=\frac{1}{n!}
\sum_{\mathcal{P}\in\mathcal{S}_n} (-1)^{|\mathcal{P}|} f(u_{\mathcal{P}(1)},\ldots,u_{\mathcal{P}(n)})$, and $\le 0$ means that one keeps only terms containing only
negative powers of the $u_i$ since otherwise they do not contribute to the
contour integral.
Now expanding the product in brackets we notice that all monomials
are of the form $\prod_i u_i^{-\beta_i}$ where $0\le\beta_i\le\alpha_i$ for all
$i$. 
\[
\mathrm{AS}\left[
\frac{\prod_{0\le i<j\le n-1}(1+u_j+u_i u_j)}{\prod_{i=0}^{n-1} u_i^{\alpha_i}}
\right]_{<0}
=\sum_{0\le\beta_i\le\alpha_i} c_{\beta}\, \det_{i,j}(u_i^{-\beta_j})
\]
where the $c_\beta$ are some coefficients.
Using the same reordering argument of the $\beta_j$ as in lemma
\ref{exptilde} (and absorbing the resulting sign in the $c_\beta$),
we conclude that
\begin{align*}
A_{\sigma,\alpha,\tau}&=\sum_{\beta\subset\alpha}c_\beta
\oint \prod_{i=0}^{n-1} \frac{d u_i}{2\pi i u_i}
\tilde s_\sigma(u) s_\tau(u) s_{\beta}(u^{-1})
\prod_{i=0}^{n-1} (1+u_i)^{n-1}
\Delta(u) \Delta(u^{-1})\\
&=\sum_{\beta\subset\alpha}c_\beta
\left(
\tilde s_\sigma(u) s_\tau(u)\prod_{i=0}^{n-1} (1+u_i)^{n-1}
\Bigg|
s_\beta(u)
\right)
\end{align*}
where by definition $(X(u)|s_\beta(u))$ is the coefficient of $s_\beta(u)$ 
in the expansion of $X(u)$ as a sum of Schur functions.
To perform this expansion, we use lemma \ref{exptilde}
and the well-known property of the multiplication of Schur functions that
the resulting Ferrers diagrams must always contain the original
ones, that is here $\sigma'$ and $\tau$. We then find that non-zero
terms are of the form $\sigma',\tau\subset\beta\subset\alpha$.
\end{proof}

Compare with lemma 4.1(1) of \cite{CKLN}, also present as lemma 3.6(a) in
\cite{Thapper}.
As a corollary,
when $\sigma\notin \A$ or $\tau\notin \A$, $A_{\sigma,\alpha,\tau}=0$.
From now on, we shall consider $A_{\sigma,\alpha,\tau}$ as a tensor
where all three indices live in $\A$.

\subsection{Connection to FPLs}
The introduction of these quantities is motivated by the
\begin{conj}\label{mainconj}
$\sum_{\alpha\in\A}
A_{\sigma,\alpha,\tau}(\mathbf{P}^{-1})_{\ \pi}^{\alpha}$ is equal to
$a_{\sigma,\pi,\tau}$,
the number of FPL configurations in a triangle with boundary conditions
specified by $\sigma$, $\pi$, $\tau$ (cf section \ref{secfpl}).
\end{conj}
This conjecture has been checked numerically up to $n=5$, using the
numerical data kindly provided by the author of \cite{Thapper}.

\subsection{Summation}
We now show how summing the $A_{\sigma,\alpha,\tau}$ according to the
prescription of \cite{Thapper} 
produces a formula which was proved in \cite{artic41} 
in the context of the $q$KZ equation, that is
on the other side of the Razumov--Stroganov conjecture.

\begin{prop} Conjecture~\ref{mainconj} implies the Razumov--Stroganov
conjecture.
\end{prop}
\begin{proof}
If conjecture~\ref{mainconj} is true, then we can compute the number
of FPL configurations with a given connectivity as follows.
Fix an arbitrary integer $k$ and define 
\begin{equation}\label{summation}
\Psi_\alpha:=\sum_{\sigma,\tau\in \A} A_{\sigma,\alpha,\tau} P_{\sigma'}(-k) P_{\tau'}(k-n+1)
\end{equation}
(compare with \eqref{presummation}).

Inserting the formula \eqref{defA} for $A_{\sigma,\alpha,\tau}$ yields
\begin{equation}
\Psi_\alpha=\braket{\prod_{i=0}^{n-1} (1+u_i)^{n-1}
\sum_{\sigma\in\A} \tilde s_\sigma(u) P_{\sigma'}(-k)
\sum_{\tau\in \A} s_\tau(u) P_{\tau'}(k-n+1)}_\alpha
\end{equation}
The summations over $\sigma$ and $\tau$ can be easily performed 
(for example by use of the dual Cauchy identity or usual Cauchy identity
depending on the sign of the arguments of $P$).
We find
$\sum_{\tau\in \A} s_\tau(u) P_{\tau'}(k-n+1)=\prod_i (1+u_i)^{k-n+1}$.
Similarly
$\sum_{\sigma\in\A} \tilde s_\sigma(u) P_{\sigma'}(-k)=
\inv(\prod_i (1+u_i)^{-k})$ but by \eqref{selfdual} we can remove $\inv$
and the formula simplifies to
\begin{equation}\label{psiformula}
\Psi_\alpha=\braket{1}_\alpha=
\Delta(u)
\prod_{0\le i<j\le n-1}(1+u_j+u_i u_j)\ \Big|_{\prod_{i=0}^{n-1} u_i^{\alpha_i}}
\end{equation}
which coincides with the expression \eqref{prepsiformula} for $\Psi'_\alpha$.

At this stage, to produce the number of FPL configurations with connectivity
$\pi$, we should then apply the matrix $\mathbf{P}^{-1}$ to the entries $\Psi_\alpha$ (by combining the
conjecture \ref{mainconj} with formula \eqref{presummation}).
However this is not necessary, since showing $\Psi'_\alpha=\Psi_\alpha$ is equivalent to $\psi'_\pi = \psi_\pi$.
\end{proof}

\subsection{Special cases}\label{specas}
The trivial case is when $\alpha=\zeron$, that is the sequence 
$(0,1,2,\ldots,n-1)$.
Since $|\alpha|=0$ we find
$A_{\sigma,\zeron,\tau}=\delta_{\sigma,\zeron}
\delta_{\tau,\zeron}$ and $\Psi_\zeron=1$.

Another special case is when $\alpha=\onen$, the largest element of $\A$,
that is the sequence $(0,2,\ldots,2n-2)$.
In this case we can use
the result first conjectured in \cite{artic41} and proved in \cite{Zeil-qKZ,artic45}:
\begin{multline*}
\left[\prod_{0\le\ell\le m\le n-1}(1-u_\ell u_m)\ 
\mathrm{AS}\left(\prod_{\ell=0}^{n-1} u_\ell^{-2\ell} \prod_{0\le\ell<m\le n-1}
(1+u_\ell u_m+\t u_m)\right)\right]_{\le 0}
\\=\mathrm{AS}\left(\prod_{\ell=0}^{n-1}
\left(u_\ell^{-1}(\t +u_\ell^{-1})\right)^\ell\right)=\prod_{0\le\ell<m\le n-1}(u_m^{-1}-u_\ell^{-1})
(\t+u_\ell^{-1}+u_m^{-1})
\end{multline*}
Here $\t=1$. We thus find
\begin{align*}
A_{\sigma,\onen,\tau}&= 
\tilde s_\sigma(u) s_\tau(u) \Delta(u)
\prod_{0\le i\le j\le n-1} (1-u_i u_j)^{-1} \prod_{i=0}^{n-1} u_i^{-2i} (1+u_i)^{n+i-1}
\big|_{CT}
\\
&=
\tilde s_\sigma(u) s_\tau(u) \Delta(u)\Delta(u^{-1})
\prod_{0\le i\le j\le n-1} (1-u_i u_j)^{-1} \prod_{0\le i<j\le n-1}(1+u_i^{-1}+u_{j}^{-1})
\prod_{i=0}^{n-1}(1+u_i)^{n-1}
\big|_{CT}
\end{align*}
(where $CT$ means extracting the constant term).
Equivalently, we have
\begin{equation*}
A_{\sigma,\onen,\tau}= \left(
\tilde s_\sigma(u) s_\tau(u)
\prod_{i=0}^{n-1} (1+u_i)^{n-1}
\prod_{0\le i\le j\le n-1} (1-u_i u_j)^{-1} 
\Bigg|
\prod_{0\le i<j\le n-1} (1+u_i+u_j)
\right)
\end{equation*}
where we have used as before the usual scalar product $(\cdot|\cdot)$
for which Schur functions are orthonormal.
Note that both bra and ket depend explicitly on $n$.

It is shown in \cite{artic41} by direct calculation 
that $\mathbf{P}^\pi_{\ \onen}=\delta^{\pi}_\onen$.
Therefore, $a_{\sigma,\onen,\tau}=A_{\sigma,\onen,\tau}$ and the formula
above gives a closed expression for it.
It seems unlikely that this expression can be simplified further.

\section{Analogues of Thapper's conjectures}
In this section we prove Thapper's conjectures
assuming conjecture \ref{mainconj}, that is we prove
some properties of the $A_{\sigma,\alpha,\tau}$ that are analogous
to those found by Thapper in the context of FPL enumeration.
\subsection{The action of $\alg$}\label{secalg}
Consider the operator that acts by inserting a factor
$s_\lambda(u)$ in the generating series. It is convenient to introduce
another bracket notation $\braket{\cdot}_{\sigma,\alpha,\tau}$
which, compared to $\braket{\cdot}_\alpha$, 
incorporates the factors $\tilde s_\sigma(u) s_\tau(u) \prod_{i=0}^{n-1}
(1+u_i)^{n-1}$.
In other words let us compute
\begin{multline}
\braket{s_\lambda(u)}_{\sigma,\alpha,\tau} := 
s_\lambda(u)\tilde s_\sigma(u) s_\tau(u) 
\Delta(u)
\prod_{i=0}^{n-1} (1+u_i)^{n-1}
\prod_{0\le i<j\le n-1}(1+u_j+u_i u_j)\ \Big|_{\prod_{i=0}^{n-1} u_i^{\alpha_i}}
\end{multline}
Similarly consider the insertion of $\tilde s_\lambda(u)$:
\begin{multline}
\braket{\tilde s_\lambda(u)}_{\sigma,\alpha,\tau} :=
\tilde s_\lambda(u)\tilde s_\sigma(u) s_\tau(u) 
\Delta(u)
\prod_{i=0}^{n-1} (1+u_i)^{n-1}
\prod_{0\le i<j\le n-1}(1+u_j+u_i u_j)\ \Big|_{\prod_{i=0}^{n-1} u_i^{\alpha_i}}
\end{multline}

Note that one can expand $s_\lambda(u)s_\tau(u)$ 
in terms of Schur functions, resulting in
\begin{equation}\label{firstid}
\braket{s_\lambda(u)}_{\sigma,\alpha,\tau}
=\sum_{\mu\in\A} \LR^\mu_{\lambda,\tau}
A_{\sigma,a,\mu}
\end{equation}
where the summation can be restricted to $\A$ because all other terms
vanish.
Of course the matrices $\LR(\lambda)_{\ \tau}^{\mu}:=
\LR^\mu_{\lambda,\tau}$ form a representation
$\lambda\mapsto\LR(\lambda)$ of $\alg$ (regular representation).

One can also expand $\tilde s_\lambda(u)s_\tau(u)$ in terms of Schur functions
$s_\mu(u)$,
with some {\it a priori}\/ unknown coefficients $\LRt^\mu_{\lambda,\tau}$:
\begin{equation}\label{secondid}
\braket{\tilde s_\lambda(u)}_{\sigma,\alpha,\tau}
=\sum_{\mu\in\A} \LRt^\mu_{\lambda,\tau}
A_{\sigma,a,\mu}
\end{equation}
One must be careful that $\LRt^\mu_{\lambda,\tau}$ is not symmetric
in the exchange of $\lambda$ and $\tau$.
If we call
$\LRt(\lambda)_{\ \tau}^{\mu}:=
\LRt^\mu_{\lambda,\tau}$,
this gives another distinct action of $\alg$. It commutes
with the previous one.

In fact in what follows it will be more natural to consider
the transpose matrices $\LR(\lambda)^T$ and 
$\LRt(\lambda)^T$,
which also form representations of $\alg$ (since it is commutative).

Of course we have dual statements by expanding this time 
the expression in brackets multiplied
by $\tilde s_\sigma(u)$ in terms of $\tilde s_\mu(u)$ and using the
fact that $\inv$ is an involution:
\begin{equation}\label{thirdid}
\braket{s_\lambda(u)}_{\sigma,\alpha,\tau}
=\sum_{\mu\in\A} \LRt^\mu_{\lambda,\sigma}
A_{\mu,a,\tau}
\end{equation}
and
\begin{equation}\label{fourthid}
\braket{\tilde s_\lambda(u)}_{\sigma,\alpha,\tau}
=\sum_{\mu\in\A} \LR^\mu_{\lambda,\sigma}
A_{\mu,a,\tau}
\end{equation}

Finally we can also expand $s_\lambda(u)\prod u_i^{-\alpha_i}$ in monomials.
This is actually more subtle than it seems because it results in monomials
for which the sequence of inverse powers is not increasing.
It is shown in appendix \ref{Kmat} (lemma \ref{Kmatex}) that
one can always reexpress the coefficient of any monomial as a linear
combination of the coefficients of the monomials with increasing inverse powers.
Note that we can discard any positive powers because they do not contribute.

That is, there exist coefficients $\mathbf{C}^\beta_{\lambda,\alpha}$ 
such that
\begin{equation}\label{existC}
\braket{s_\lambda(u)F(u)}_\alpha=\sum_{\beta\in\A} 
\mathbf{C}_{\lambda,\alpha}^\beta
\braket{F(u)}_\beta
\end{equation}
for any symmetric function $F$, 
and similarly for $\mathbf{\tilde C}^\beta_{\lambda,\alpha}$.

In the present case, we find
\begin{equation}\label{Crel}
\braket{s_\lambda(u)}_{\sigma,\alpha,\tau}
=\sum_{\beta\in\A} \mathbf{C}^\beta_{\lambda,\alpha}
A_{\sigma,\beta,\tau}
\end{equation}
and
\begin{equation}\label{sixthid}
\braket{\tilde s_\lambda(u)}_{\sigma,\alpha,\tau}
=\sum_{b\in\A} \mathbf{\tilde C}^\beta_{\lambda,\alpha}
A_{\sigma,\beta,\tau}
\end{equation}
which can be made into matrices
$\mathbf{C}(\lambda)_{\ \alpha}^{\beta}:=\mathbf{C}^\beta_{\lambda,\alpha}$ 
and $\mathbf{\tilde C}(\lambda)_{\ \alpha}^{\beta}:=\mathbf{\tilde C}^\beta_{\lambda,\alpha}$.

Combining these relations (\ref{firstid}--\ref{sixthid})
results in various identities.

A special case occurs when one considers the multiplication
by $\prod_i(1+u_i)$, because of the fact that it is invariant by
the involution. We have:
\begin{equation}\label{mainthapper}
\braket{\prod_i (1+u_i)}_{\sigma,\alpha,\tau}
=\sum_{\mu\in\A} \LR^\mu_{\ \tau} A_{\sigma,\alpha,\mu}
=\sum_{\mu\in\A} \LR^\mu_{\ \sigma} A_{\mu,\alpha,\tau}
=\sum_{\beta\in\A} \mathbf{C}^\beta_{\ \alpha} A_{\sigma,\beta,\tau}
\end{equation}
with the obvious notations 
$\LR^\mu_{\ \tau}=\sum_{i=0}^{n-1}\LR^\mu_{e_i,\tau}$ 
and
$\mathbf{C}^\beta_{\ \alpha}=\sum_{i=0}^{n-1}\mathbf{C}^\beta_{e_i,\alpha}$.

Observe also that taking $\sigma=\zeron$ in \eqref{fourthid} and 
combining with \eqref{secondid} results in
\begin{equation}\label{otherthapper}
A_{\lambda,\alpha,\tau}
=\sum_{\mu\in\A} \LRt^\mu_{\lambda,\tau} A_{\zeron,\alpha,\mu}
\end{equation}

\subsection{Matrix identities}
Introduce, following \cite{Thapper}, the matrices $A(\lambda)_{\tau \alpha}:=A_{\lambda,\alpha,\tau}$
and $\bar A(\alpha)_{\sigma\tau}=A_{\sigma,\alpha,\tau}$.
Note that lemma~\ref{uptri} says that the $A(\sigma)$ are upper triangular.
We now rewrite the various identities of the previous section in these
matrix notations.

Relations (\ref{firstid}--\ref{sixthid}) become
\begin{align}\label{firstidm}
\LRt(\lambda)^T \bar A(\alpha)&=
\bar A(\alpha) \LR(\lambda)\\\label{secondidm}
\LR(\lambda)^T A(\sigma)&=A(\sigma)\mathbf{C}(\lambda)\\
\LR(\lambda)^T \bar A(\alpha)&=
\bar A(\alpha) \LRt(\lambda)\\
\LRt(\lambda)^T A(\sigma)&=A(\sigma)\mathbf{\tilde C}(\lambda)
\label{lastidm}\end{align}
and \eqref{mainthapper} becomes
\begin{align}
\LR^T \bar A(\alpha)&=\bar A(\alpha)\LR\\
\LR^T A(\alpha)&=A(\alpha)\mathbf{C}\label{nadeauthm}
\end{align}
Compare with conj.~3.4 and 3.5(a) of \cite{Thapper} (with $B=\LR^T$:
$B$ removes columns of boxes to Ferrers diagrams whereas $\LR$
adds them -- so that $B$ is upper triangular whereas $\LR$
is lower triangular).

As to \eqref{otherthapper}, it becomes
\begin{equation}
A(\lambda)=\LRt(\lambda)^T A(\zeron)
\end{equation}
According to lemmas \ref{uptri} and \ref{nadeauthmb}, $A(\zeron)$ is upper
triangular with ones on the diagonal, and therefore invertible.
The matrices $A(\lambda)A(\zeron)^{-1}$ are equal to
$\LRt(\lambda)^T$
and thus satisfy the relations of the $\alg$ algebra, 
cf prop.~3.9 of \cite{Thapper}. In particular
they commute, cf prop.~3.10 of \cite{Thapper}.

One can similarly combine \eqref{fourthid} and \eqref{sixthid} to obtain
\begin{equation}
A(\lambda)=A(\zeron)\mathbf{\tilde C}(\lambda)
\end{equation}
i.e.\ $\mathbf{\tilde C}(\lambda)=A(\zeron)^{-1}A(\lambda)$.

More directly, we conclude from (\ref{secondidm},\ref{lastidm}) that 
the matrices $\mathbf{C}(\lambda)$ and $\mathbf{\tilde C}(\lambda)$
provide us with another pair of commuting representations of $\alg$; and that
$A(\zeron)$ intertwines the representations
$\mathbf{C}(\lambda)$ and $\LR(\lambda)^T$,
as well as
$\mathbf{\tilde C}(\lambda)$ and $\LRt(\lambda)^T$.
In fact, any linear combination of the $A(\sigma)$ is an intertwiner,
but for it to be invertible the coefficient of $A(\zeron)$ must be non-zero.

\subsection{The involution}
Consider the matrix $\matinv$ of the involution $\inv$ acting on $\alg$:
\begin{equation}
\tilde s_\lambda=\sum_{\mu\in\A} \matinv_{\ \lambda}^{\mu} s_\mu+\text{terms not
inside $\onen$}
\end{equation}
(equivalently $\matinv_{\ \lambda}^{\mu}:=\LRt^\mu_{\lambda,\zeron}$).
Then by definition we have
\begin{equation}
\LRt(\lambda)=
\matinv
\LR(\lambda)
\matinv
\end{equation}
That is, $\matinv$ intertwines the representations 
$\LR(\lambda)$
and
$\LRt(\lambda)$
of $\alg$.

We now have a chain of (invertible) intertwiners between the four representations
$\mathbf{C}(\lambda)$,
$\LR(\lambda)^T$, $\LRt(\lambda)^T$,
$\mathbf{\tilde C}(\lambda)$.
Composing them, we find that
\begin{equation}
\mathbf{\tilde C}(\lambda)=\mathbf{R} \mathbf{C}(\lambda)\mathbf{R}
\end{equation}
with $\mathbf{R}:=A(\zeron)^{-1}\matinv^T A(\zeron)$, $\mathbf{R}^2=1$.

What is the meaning of $\mathbf{R}$? It is a simple exercise to check, combining
the various relations found so far, that
\begin{equation}
\sum_{\beta\in\A} 
A_{\sigma,\beta,\tau}
\mathbf{R}_{\ \alpha}^\beta 
=A_{\tau,\alpha,\sigma}
\end{equation}
\rem{find this explicitly just using the CT definition}
That is, $\mathbf{R}$ corresponds to the operation of mirror image in the space
of link patterns. In the canonical link pattern basis ($\mathbf{r}:=\mathbf{P} \mathbf{R}\mathbf{P}^{-1}$), we simply
have $\mathbf{r}^{\pi}_{\ \pi'}=1$ if $\pi$ and $\pi'$ are mirror images of each other,
and $0$ otherwise.

\subsection{Recurrence}\label{secrecur}
\newcommand\surr[1]{{(#1)_m}}
Let $m$ be a non-negative integer.
By definition, let $\surr{\alpha}$ be the sequence of increasing integers in $\Am$
associated to a sequence $\alpha\in\A$ by the embedding
of $\A$ into $\Am$ described in section \ref{secembed}. Explicitly,
$\surr{\alpha}=(0,1,\ldots,m-1,m+\alpha_0,\ldots,m+\alpha_{n-1})$.
In other words, viewed as a Dyck path, $\surr{\alpha}$ consists of $m$ up steps,
then the Dyck path $\alpha$ of length $2n$, then $m$ down steps. 
We also write $(\alpha)_1=(\alpha)$.
Let us derive a recurrence relation for $\Psi_\surr{\alpha}$ 
using the formalism of the previous sections.

Start with the formula \eqref{psiformula} for $\Psi_\surr{\alpha}$ and note that
we are simply looking for the constant term
as a function of $u_0$. So we can set $u_0=0$; as a result, the Vandermonde
determinant shifts all powers of $u_i$ by 1, so that one can now apply
the same argument to $u_1$, etc.
After $m$ steps we find the following formula:
\begin{equation}
\Psi_{\surr{\alpha}}=
\Delta(u)
\prod_{i=0}^{n-1} (1+u_i)^{m}
\prod_{0\le i<j\le n-1}(1+u_j+u_i u_j)\ \Big|_{\prod_{i=0}^{n-1} u_i^{\alpha_i}}
\end{equation}

Using bracket notations, this is nothing but
\begin{equation}
\Psi_{\surr{\alpha}}=
\braket{\prod_{i=0}^{n-1} (1+u_i)^{m}}_\alpha
\end{equation}
Note that this means that $\Psi$
does {\em not}\/ satisfy a stability property with respect to the embedding
of $\A$ into $\Am$. 

The insertion of such products has already been analyzed in 
section \ref{secalg}.
Applying \eqref{existC}, we find
a recurrence relation for $\Psi_\alpha$:
\begin{equation}\label{recura}
\Psi_{\surr{\alpha}}=\sum_{\beta\in \A} \Psi_\beta (\mathbf{C}^m)^\beta_{\ \alpha} 
\end{equation}
It formulates a component of size $m+n$ in terms of components of size $n$.

Let us now discuss briefly the change of basis to the link pattern
basis. It turns out to be ``compatible'' with the operation
$\surr{\cdot}$ in the following sense: 
for all $\rho\in \Am$, $\alpha\in\A$, 
\begin{equation}
\mathbf{P}_{\ \surr{\alpha}}^{\rho}=
\begin{cases}
\mathbf{P}_{\ \alpha}^{\pi}&\text{if $\rho=\surr{\pi}$ for some $\pi\in \A$}\\
0&\text{otherwise}
\end{cases}
\end{equation}
(this property is a combination of the upper triangularity and stability
with respect to the embedding of $\mathbf{P}$).
Thus, writing
\begin{equation}
\Psi_{\surr{\alpha}}=\sum_{\pi\in \A} \Psi_{\surr{\pi}} \mathbf{P}^{\pi}_{\ \alpha} 
\end{equation}
and inverting $\mathbf{P}$,
we can rewrite \eqref{recura}
\begin{equation}\label{recur}
\psi_{\surr{\pi}}=\sum_{\rho\in \A} \psi_\rho (\mathbf{c}^m)^\rho_{\ \pi} 
\end{equation}
with $\mathbf{c}:=\mathbf{P}\mathbf{C}\mathbf{P}^{-1}$.

In particular, for $m=1$, since any link pattern
contains a pairing of neighbors, which after appropriate rotation can
be mapped to $(0,2n-1)$, \eqref{recur} provides a closed recurrence relation
for the $\psi_\pi$.

NB: one can also write recurrence relations for the $A_{\sigma,\alpha,\tau}$.
Following the same reasoning as for $\Psi$, but paying
attention to the extra factors occurring because of the explicit
dependence of $A_{\sigma,\surr{\alpha},\tau}$ on the size $m+n$, we find that
$A_{\sigma,\surr{\alpha},\tau}
=\braket{\prod_{i=0}^{n-1}(1+u_i)^{2m}}_{\sigma,\alpha,\tau}$.
Just like $\Psi$, $A_{\sigma,\alpha,\tau}$ does not
satisfy a stability property with respect to the embedding
of $\A$ into $\Am$. Use of relations \eqref{mainthapper} give various
identities, including recurrences, for the $A_{\sigma,\alpha,\tau}$.
Also, if one uses
one of the first two equalities of \eqref{mainthapper} and apply it to
the relation \eqref{summation}, one finds
\begin{equation*}
\Psi_{\surr{\alpha}}=\sum_{\sigma,\tau\in \A} A_{\sigma,\alpha,\tau}
P_{\sigma'}(2m-k) P_{\tau'}(k-n-m+1)
\end{equation*}
or the same expression with instead $P_{\sigma'}(-k)P_{\tau'}(k-n+m+1)$,
which is of course the same up to a change of the arbitrary
integer $k$. This expression is Eq.~(4) (with arbitrary $m$) 
of \cite{Thapper} (see also Theorem 4.2 of \cite{CKLN}).

\subsection{Largest component}
As has already been mentioned in section \ref{specas},
the change of basis is trivial for the largest component, so
that the component $\psi_\onen$
(i.e.\ with $\onen$ viewed as a link
pattern) is the same as the component $\Psi_\onen$ in our basis (i.e.\ 
the component with sequence $\alpha_i=2i$). 

Furthermore, it is also known \cite{artic41}
 how to obtain the sum of components;
one has
\begin{equation*}
\sum_{\epsilon_1,\ldots,\epsilon_{n-1}\in\{0,1\}}
\mathbf{P}^{\pi}_{\ (0,2-\epsilon_1,\ldots,2(n-1)-\epsilon_{n-1})}
=1\qquad\forall \pi\in \A
\end{equation*}
Note that this is exactly the image of $\onen$ by the action of $\mathbf{C}$,
so that $(\mathbf{P}\mathbf{C})^\pi_{\ \onen}=1$ for all $\pi$. 
Multiplying on the right by $\mathbf{P}^{-1}$ and using again
$\mathbf{P}^{\pi}_{\ \onen}=\delta^{\pi}_\onen$ we find
$\mathbf{c}^\pi_{\ \onen}=1$,
which is exactly the content of conj.~3.5(b) of \cite{Thapper}.

Thus, as explained in \cite{Thapper}, the recurrence produces in this case the identity
\begin{equation}
\psi_{\onenp}=\sum_{\pi\in \A}\psi_\pi
\end{equation}
as well as the more general one
\begin{equation}
\psi_{(\onen)_{m+1}}=\sum_{\pi\in \A}\psi_{(\pi)_m}
\end{equation}
to compare with the conjectures of \cite{Zuber-conj}.

\subsection{The $A$'s from the $\Psi$'s}
It it perhaps suggestive that using the various matrices defined above,
one can in fact compute $A_{\sigma,\alpha,\tau}$ from the data of the $\Psi_\alpha$
alone. Indeed, introducing yet another notation
$A(\sigma,\tau)$ for the row vector with entries $A_{\sigma,\alpha,\tau}$,
and similarly for $\Psi$, we have
\begin{equation}
A(\sigma,\tau)=\Psi\,\mathbf{\tilde C}(\sigma)\mathbf{C}(\tau)
\mathbf{C}^{n-1} 
\end{equation}

\section{The $\t$-generalization}
In \cite{Pas-RS,artic34}, it was suggested how to generalize the ground state
eigenvector $\psi'$ of the Temperley--Lieb(1) model into a vector
depending on an extra parameter $\t$ (often written as
$\t=-q-q^{-1}$) and which is obtained
by specializing a certain polynomial solution of the 
quantum Knizhnik--Zamolodchikov ($q$KZ) 
equation. The original vector is recovered by
taking $\t=1$. The parameter $\t$ does not seem to
have any obvious meaning in terms of FPLs, though a connection
to Totally Symmetric Self-Complementary Plane Partitions was found in \cite{DF-qKZ-TSSCPP}.
Many formulae generalize to $\t$ away from unity \cite{artic41}.
It is natural to wonder if the formulae of the present work fall in this
category. 
We briefly present here this generalization.

Introduce the $\t$-dependent bracket
\[
\braket{F(u)}_\alpha:=
F(u)\Delta(u)
\prod_{0\le i<j\le n-1}(1+\t u_j+u_i u_j)\ \Big|_{\prod_{i=0}^{n-1} u_i^{\alpha_i}}
\]
This produces a natural generalization of the $\Psi_a$: 
\[
\Psi_a = \braket{1}_\alpha
\]
which corresponds to the
homogeneous limit of a solution of the $q$KZ system of difference equations.
The change of basis to the link pattern basis can be found for all
$\t$ in \cite{artic41}
(or in \eqref{formulaP} where the $U_i$ are Chebyshev polynomials of $\t$).

Since the interpretation in terms of FPLs in a
triangle is lacking, there seems little point in introducing the
$\t$-generalization of $A_{\sigma,a,\tau}$. However, one can still define the
various matrices of multiplication by a Schur function $s_\lambda$.
The matrices $\LR(\lambda)$ are of course unchanged; however the
matrices $\mathbf{C}(\lambda)$ which act on increasing sequences
by 
\[
\braket{s_\lambda(u)F(u)}_\alpha=\sum_{\beta\in\A} 
\mathbf{C}_{\lambda,\alpha}^\beta
\braket{F(u)}_\beta
\]
now have entries which are polynomials in $\t$,
and can be computed using the methods
of appendices \ref{Kmat} and 
\ref{Pmat}. They still form a reprentation of $\alg$.
The dual versions can be defined via the involution
\[
\inv\left(\prod_i (1+z u_i)\right)=\prod_i \frac{1}{1-z\frac{u_i}{1+\t u_i}}
\]

Of particular interest is the (self-dual) operator $\mathbf{C}$
of multiplication by $\prod_{i=0}^{n-1}(1+\t u_i)$, that is
\[
\mathbf{C}=\sum_{i=0}^{n-1} \mathbf{C}(e_i)\t^i
\]
It corresponds to setting $z=\t$ in the notations of appendix \ref{Pmat}.

Then all the recursion relations conjectured in \cite{Thapper} are
satisfied by $\Psi$ for all $\t$, with the use of this modified 
matrix $\mathbf{C}$.
Namely, equations \eqref{recura} and \eqref{recur} hold without any change.
However, naively, they do not form a closed set of recursion relations.
The reason is that for $\t\ne 1$, Wieland's rotational invariance
theorem does not hold any more: in general $\psi_\pi \ne \psi_{\rho(\pi)}$,
so one cannot assume that there exists a pairing $(0,2n-1)$.

However, a remarkable phenomenon occurs:
in the link pattern basis, one has
\[
\psi_{\pi}=\sum_{\rho\in \A} \psi_\rho \mathbf{c}^{(\rho)}_{\ \ \pi} 
\]
where now $\pi$ is an arbitrary link pattern of size $n+1$, and $\mathbf{c}$
is the matrix of size $c_{n+1}$. This is a non-trivial generalization of
Eq.~\eqref{recur} at $m=1$, in the sense that in the special
case where $\pi$ has a pairing $(0,2n-1)$, i.e.\ $\pi=(\pi')$,
we recover \eqref{recur}.

Thus, we have a closed recursion again,
and the vectors $\psi$ for successive values of $n$ can be obtained
by simply iterating the matrices $\mathbf{c}$ in size $1,\ldots,n$.
This observation will be the subject of future work.

\bigskip
\noindent {\it Note added.}
After this work was completed,
Nadeau announced
\cite{Nadeau-FPL1}
a direct proof of conj.~3.4 of \cite{Thapper}
i.e.\ the analogue of \eqref{nadeauthm}. He also announced
\cite{Nadeau-FPL2}
a bijection between Knutson--Tao puzzles
and FPLs in a triangle with $|\alpha|=|\sigma|+|\tau|$, which would
provide a bijective proof of the second part of lemma \ref{nadeauthmb}.

\appendix

\section{Existence of the matrix $\mathbf{K}$}\label{Kmat}
In this section we no longer assume that the
sequences $(\alpha_i)$ are increasing.
\begin{lemma}\label{Kmatex}
Let $\alpha=(\alpha_0,\ldots,\alpha_{n-1})$ be an arbitrary sequence of integers
such that $\alpha_i\le 2i$. There exist coefficients $\mathbf{K}_{\ \alpha}^\beta$ such that
for any symmetric function $F(u)$,
\begin{equation}\label{lemmarel}
\braket{F(u)}_\alpha=\sum_{\beta\in \A} \braket{F(u)}_\beta \mathbf{K}_{\ \alpha}^\beta 
\end{equation}
\end{lemma}

\begin{proof}
The relation can be checked on Schur functions $s_\tau(u)$ only.
Furthermore, if $\tau\notin \A$, because $\alpha_i\le 2i$,
both sides of the equality are zero. So one can think of $\braket{\cdot}_\alpha$
as a linear form on $\alg$.
It it actually convenient at this stage to switch to a different basis of
$\alg$, namely $s_\tau(u) \prod_{i=0}^{n-1}(1+u_i)^{n-1}$ (in order to avoid
having to introduce yet another new notation).
Plugging this into \eqref{lemmarel}, we note that
the right hand side is simply $(A(\zeron)\mathbf{K})_{\tau,\alpha}$.
Let us thus define
\begin{equation}
A^{ext}(\zeron)_{\tau,\alpha}=\braket{s_\tau(u)\prod_{i=0}^{n-1}(1+u_i)^{n-1}}_\alpha
\end{equation}
that is the same definition as for $A$, but in which we relax the
constraint that $(\alpha_i)$ is increasing. Then, by definition,
\begin{equation}
\mathbf{K}=A(\zeron)^{-1}A^{ext}(\zeron)
\end{equation}
satisfies the relation \eqref{lemmarel}.
\end{proof}

This way, we see that we can build
a matrix $\mathbf{C}(\lambda)$ that satisfies \eqref{existC}: simply define
\begin{equation}\label{Cext}
\mathbf{C}(\lambda)=\mathbf{K} \mathbf{C}^{ext}(\lambda)
\end{equation}
where $\mathbf{C}^{ext}(\lambda)$ is an explicit matrix which encodes
the decomposition of $s_\lambda(u)$ into monomials; for example
in the case of $\mathbf{C}^{ext}_z:=\sum_{i=0}^{n-1}\mathbf{C}^{ext}(e_i)z^i$
corresponding to multiplication by $\prod_i (1+z u_i)$, we get
\begin{equation}\label{Cextz}
\mathbf{C}^{ext}_z\, ^\beta_{\ \alpha}=
\begin{cases}
z^{\sum_i(\alpha_i-\beta_i)}&\text{if $\alpha_i-\beta_i\in\{0,1\}$ for all $i$}\\
0&\text{otherwise}
\end{cases}
\end{equation}
The same procedure works of course for $\mathbf{\tilde C}(\lambda)$.

This construction of $\mathbf{C}(\lambda)$, $\mathbf{\tilde C}(\lambda)$
is not very explicit because it assumes the knowledge of $\mathbf{K}$,
which requires to compute $A^{ext}(\zeron)$
and to invert $A(\zeron)$. We provide a simpler formula in the
next section.

\section{The matrices $\mathbf{P}$, $\mathbf{P}^{ext}$}\label{Pmat}
We describe the change of basis from increasing sequences to link patterns,
given by the matrix $\mathbf{P}$. In fact we describe a slightly bigger
matrix, $\mathbf{P}^{ext}$, which allows for more general 
(non-decreasing) sequences.
This is borrowed from appendix A of
\cite{artic41} and since these results are not needed in the rest of the paper,
we do not provide their proof and refer to \cite{artic41} for details.

Define $U_i:=1,-1,0$ depending on whether
$i=0,1,2$ mod 3 (these are evaluations of Chebyshev polynomials at
the parameter $\t=1$). Given a sequence of integers $\alpha=(\alpha_i)$, define
\begin{equation}\label{formulaP}
\mathbf{P}^\pi_{\ \alpha}=
\prod_{\scriptstyle 0\le i<j<2n\atop\scriptstyle
\text{$i$ and $j$ paired by $\pi$}} U_{\# \{ \ell: i\le \alpha_\ell<j\} - (j-i+1)/2}
\end{equation}
Let us call $\mathbf{P}$ the matrix of $\mathbf{P}^\pi_{\ \alpha}$ 
in which $\alpha\in A_n$, and
$\mathbf{P}^{ext}$ to be given by
\begin{equation}\label{Pext}
\mathbf{P}^{ext}=
\mathbf{P}\mathbf{K}
\end{equation}
that is the matrix of passage from {\em any}\/ sequence $\alpha$ to link
patterns $\pi$ (see appendix \ref{Kmat}).
The claim is that \eqref{formulaP} provides the entries of $\mathbf{P}^{ext}$
for any {\em non-decreasing}\/ sequence $\alpha$.

As an application, let us provide a direct way of computing
$\mathbf{c}(\lambda):=\mathbf{P}\mathbf{C}(\lambda)\mathbf{P}^{-1}$, 
the matrix that implements multiplication
by $s_\lambda(u)$ in the link pattern basis. Using \eqref{Cext} and
\eqref{Pext}, we find
\begin{equation}\label{Cb}
\mathbf{c}(\lambda):=\mathbf{P} \mathbf{C}(\lambda)\mathbf{P}^{-1}=\mathbf{P}^{ext}\mathbf{C}^{ext}\mathbf{P}^{-1}
\end{equation}
Consider now the generating function of elementary symmetric functions
corresponding to multiplying by $\prod_i (1+z u_i)$:
\begin{equation}
\mathbf{c}_z:=
\sum_{i=0}^{n-1}\mathbf{c}(e_i) z^i
\end{equation}
It is easy to check that starting from a monomial
$\prod_i u_i^{-\alpha_i}$ with an increasing sequence $\alpha$ and multiplying
by $\prod_i (1+z u_i)$, one produces sums of monomials with non-decreasing
sequences of powers. Thus one can use \eqref{formulaP} 
(and \eqref{Cextz}) to calculate
explicitly $\mathbf{c}_z$ using \eqref{Cb}. Any 
$\mathbf{c}(\lambda)$ can then be computed by applying say
the Von N\"agelsbach--Kostka identity (dual Jacobi--Trudi identity).
If one is only interested in $\mathbf{c}=\mathbf{c}_1$,
the matrix that appears in recurrence formulae (cf \eqref{recur})
one can of course directly set $z=1$.

\section{Example of ground state entry of the Temperley--Lieb loop model}\label{PFTL}
\smallboxes
\newcommand\ya{\tableau{\ }}
\newcommand\yb{\tableau{\ \\\ }}
\newcommand\yc{\tableau{\ &\ }}
\newcommand\yd{\tableau{\ &\ \\\ }}
We consider the model in size $2n=6$. 
The basis is ordered as follows:
$\zeron$, $\ya$, $\yb$, $\yc$,
$\onen=\yd$ (modulo the bijection to link patterns of section \ref{secbij}).
The matrices of the $\mathbf{e}_i$ are:
\begin{align*}
\mathbf{e}_1&=\begin{pmatrix}
0 & 0 & 0 & 0 & 0 \\
 0 & 0 & 0 & 0 & 0 \\
 0 & 1 & 1 & 0 & 0 \\
 0 & 0 & 0 & 0 & 0 \\
 1 & 0 & 0 & 1 & 1
\end{pmatrix}
&\qquad
\mathbf{e}_2&=\begin{pmatrix}
 0 & 0 & 0 & 0 & 0 \\
 1 & 1 & 1 & 0 & 0 \\
 0 & 0 & 0 & 0 & 0 \\
 0 & 0 & 0 & 1 & 1 \\
 0 & 0 & 0 & 0 & 0
\end{pmatrix}
\\
\mathbf{e}_3&=\begin{pmatrix}
 1 & 1 & 0 & 0 & 0 \\
 0 & 0 & 0 & 0 & 0 \\
 0 & 0 & 0 & 0 & 0 \\
 0 & 0 & 0 & 0 & 0 \\
 0 & 0 & 1 & 1 & 1
\end{pmatrix}
&\qquad
\mathbf{e}_4&=\begin{pmatrix}
 0 & 0 & 0 & 0 & 0 \\
 1 & 1 & 0 & 1 & 0 \\
 0 & 0 & 1 & 0 & 1 \\
 0 & 0 & 0 & 0 & 0 \\
 0 & 0 & 0 & 0 & 0
\end{pmatrix}
\\
\mathbf{e}_5&=\begin{pmatrix}
 0 & 0 & 0 & 0 & 0 \\
 0 & 0 & 0 & 0 & 0 \\
 0 & 0 & 0 & 0 & 0 \\
 0 & 1 & 0 & 1 & 0 \\
 1 & 0 & 1 & 0 & 1
\end{pmatrix}
&\qquad
\mathbf{e}_6&=\begin{pmatrix}
 1 & 0 & 0 & 0 & 1 \\
 0 & 1 & 1 & 1 & 0 \\
 0 & 0 & 0 & 0 & 0 \\
 0 & 0 & 0 & 0 & 0 \\
 0 & 0 & 0 & 0 & 0
\end{pmatrix}
\end{align*}
so that
\[
\mathbf{H}=\begin{pmatrix}
 2 & 1 & 0 & 0 & 1 \\
 2 & 3 & 2 & 2 & 0 \\
 0 & 1 & 2 & 0 & 1 \\
 0 & 1 & 0 & 2 & 1 \\
 2 & 0 & 2 & 2 & 3
\end{pmatrix}
\]
with Perron--Frobenius eigenvector
\[
\psi'=\begin{pmatrix} 1 & 2 & 1 & 1 & 2 \end{pmatrix}
\]

\section{Example of matrices $\bar A$ and $\mathbf{C}$}
We again provide data for $n=3$. 
The basis is ordered as in the previous section, but
in this section the matrices are as defined in sections 4 and 5 i.e.\ 
we do not use the link pattern basis.
To recover the data of \cite{Thapper} or of the previous section,
one needs to use the change of basis 
(cf \eqref{chgbasis}) given by
the matrix
\begin{equation*}
\mathbf{P}=\begin{pmatrix}
1 & 0 & 1 & 0 & 0 \\
 0 & 1 & 0 & 0 & 0 \\
 0 & 0 & 1 & 0 & 0 \\
 0 & 0 & 0 & 1 & 0 \\
 0 & 0 & 0 & 0 & 1
\end{pmatrix}
\end{equation*}

First we give the $A_{\sigma,\alpha,\tau}$ under the form of the matrices $\bar A$:
\begin{align*}
\bar A(\zeron)&=
\begin{pmatrix}
 1 & 0 & 0 & 0 & 0 \\
 0 & 0 & 0 & 0 & 0 \\
 0 & 0 & 0 & 0 & 0 \\
 0 & 0 & 0 & 0 & 0 \\
 0 & 0 & 0 & 0 & 0
\end{pmatrix}&\qquad
\bar A(\ya)&=
\begin{pmatrix}
 4 & 1 & 0 & 0 & 0 \\
 1 & 0 & 0 & 0 & 0 \\
 0 & 0 & 0 & 0 & 0 \\
 0 & 0 & 0 & 0 & 0 \\
 0 & 0 & 0 & 0 & 0
\end{pmatrix}\\
\bar A(\yb)&=
\begin{pmatrix}
 7 & 3 & 1 & 0 & 0 \\
 4 & 1 & 0 & 0 & 0 \\
 0 & 0 & 0 & 0 & 0 \\
 1 & 0 & 0 & 0 & 0 \\
 0 & 0 & 0 & 0 & 0
\end{pmatrix}&\qquad
\bar A(\yc)&=
\begin{pmatrix}
 6 & 4 & 0 & 1 & 0 \\
 3 & 1 & 0 & 0 & 0 \\
 1 & 0 & 0 & 0 & 0 \\
 0 & 0 & 0 & 0 & 0 \\
 0 & 0 & 0 & 0 & 0
\end{pmatrix}\\
\bar A(\yd)&=
\begin{pmatrix}
17 & 13 & 4 & 3 & 1 \\
 13 & 7 & 1 & 1 & 0 \\
 4 & 1 & 0 & 0 & 0 \\
 3 & 1 & 0 & 0 & 0 \\
 1 & 0 & 0 & 0 & 0
\end{pmatrix}\\
\end{align*}
Summing these objects, cf \eqref{summation}, reproduces $\Psi$:
\begin{equation*}
\Psi=\begin{pmatrix} 1 & 2 & 2 & 1 & 2 \end{pmatrix}
\end{equation*}
Note that one possible choice in \eqref{summation} is to set $k=0$;
since $P_\lambda(0)=\delta_{\lambda}^\zeron$, one sums
in this case the first lines of the matrices $\bar A$ only,
with coefficients $P_{\tau'}(-n+1)$ which is nothing but
$(-1)^{|\tau|}$ times 
the dimension of $\tau$ viewed as a $GL(n-1)$ representation.
In the present case, we find $1,-2,1,3,-2$. Of course the same linear
combination works when summing the first rows only (set $k=n-1$).

Next, let us describe the various matrices of the representations of $\alg$.
The $\LR$ are Littlewood--Richardson coefficients, so they are known.
Using the matrix
\begin{equation*}
A(\zeron)=\begin{pmatrix}
1 & 4 & 7 & 6 & 17\\
0 & 1 & 3 & 4 & 13\\
0 & 0 & 1 & 0 & 4 \\
0 & 0 & 0 & 1 & 3 \\
0 & 0 & 0 & 0 & 1 \\
\end{pmatrix}
\end{equation*}
one can build the $\mathbf{C}$:
\begin{align*}
\LR(\zeron)&=
\begin{pmatrix}
 1 & 0 & 0 & 0 & 0 \\
 0 & 1 & 0 & 0 & 0 \\
 0 & 0 & 1 & 0 & 0 \\
 0 & 0 & 0 & 1 & 0 \\
 0 & 0 & 0 & 0 & 1
\end{pmatrix}
&\qquad
\mathbf{C}(\zeron)&=
\begin{pmatrix}
1 & 0 & 0 & 0 & 0 \\
 0 & 1 & 0 & 0 & 0 \\
 0 & 0 & 1 & 0 & 0 \\
 0 & 0 & 0 & 1 & 0 \\
 0 & 0 & 0 & 0 & 1
\end{pmatrix}
\\
\LR(\ya)&=
\begin{pmatrix}
 0 & 0 & 0 & 0 & 0 \\
 1 & 0 & 0 & 0 & 0 \\
 0 & 1 & 0 & 0 & 0 \\
 0 & 1 & 0 & 0 & 0 \\
 0 & 0 & 1 & 1 & 0
\end{pmatrix}
&\qquad
\mathbf{C}(\ya)&=
\begin{pmatrix}
 0 & 1 & - 1 & 0 & 0 \\
 0 & 0 & 1 & 1 & 0 \\
 0 & 0 & 0 & 0 & 1 \\
 0 & 0 & 0 & 0 & 1 \\
 0 & 0 & 0 & 0 & 0
\end{pmatrix}
\\
\LR(\yb)&=
\begin{pmatrix}
 0 & 0 & 0 & 0 & 0 \\
 0 & 0 & 0 & 0 & 0 \\
 1 & 0 & 0 & 0 & 0 \\
 0 & 0 & 0 & 0 & 0 \\
 0 & 1 & 0 & 0 & 0
\end{pmatrix}
&\qquad
\mathbf{C}(\yb)&=
\begin{pmatrix}
0 & 0 & 1 & 0 & 0 \\
 0 & 0 & 0 & 0 & 1 \\
 0 & 0 & 0 & 0 & 0 \\
 0 & 0 & 0 & 0 & 0 \\
 0 & 0 & 0 & 0 & 0
\end{pmatrix}
\\
\LR(\yc)&=
\begin{pmatrix}
 0 & 0 & 0 & 0 & 0 \\
 0 & 0 & 0 & 0 & 0 \\
 0 & 0 & 0 & 0 & 0 \\
 1 & 0 & 0 & 0 & 0 \\
 0 & 1 & 0 & 0 & 0
\end{pmatrix}
&\qquad
\mathbf{C}(\yc)&=
\begin{pmatrix}
0 & 0 & 0 & 1 & - 1 \\
 0 & 0 & 0 & 0 & 1 \\
 0 & 0 & 0 & 0 & 0 \\
 0 & 0 & 0 & 0 & 0 \\
 0 & 0 & 0 & 0 & 0
\end{pmatrix}
\\
\LR(\yd)&=
\begin{pmatrix}
 0 & 0 & 0 & 0 & 0 \\
 0 & 0 & 0 & 0 & 0 \\
 0 & 0 & 0 & 0 & 0 \\
 0 & 0 & 0 & 0 & 0 \\
 1 & 0 & 0 & 0 & 0
\end{pmatrix}
&\qquad
\mathbf{C}(\yd)&=
\begin{pmatrix}
 0 & 0 & 0 & 0 & 1 \\
 0 & 0 & 0 & 0 & 0 \\
 0 & 0 & 0 & 0 & 0 \\
 0 & 0 & 0 & 0 & 0 \\
 0 & 0 & 0 & 0 & 0
\end{pmatrix}
\end{align*}
so that finally,
\[
\mathbf{C}=\mathbf{C}(\zeron)+\mathbf{C}(\ya)+\mathbf{C}(\yb)
=
\begin{pmatrix}
 1 & 1 & 0 & 0 & 0 \\
 0 & 1 & 1 & 1 & 1 \\
 0 & 0 & 1 & 0 & 1 \\
 0 & 0 & 0 & 1 & 1 \\
 0 & 0 & 0 & 0 & 1
\end{pmatrix}
\]

The dual versions $\LRt$, $\mathbf{\tilde C}$ can be obtained by use
of the involution whose matrix is
\begin{equation*}
\matinv=\begin{pmatrix}
1 & 0 & 0 & 0 & 0 \\
 0 & 1 & 0 & 0 & 0 \\
 0 & 1 & 0 & 1 & 0 \\
 0 & -1 & 1 & 0 & 0 \\
 0 & -1 & 1 & -1 & 1
\end{pmatrix}
\end{equation*}
and will not be listed here.

\renewcommand\MR[1]{\relax\ifhmode\unskip\spacefactor3000 \space\fi
  \MRhref{#1}{{\sc mr}}}
\renewcommand{\MRhref}[2]{%
  \href{http://www.ams.org/mathscinet-getitem?mr=#1}{#2}}
\bibliography{../biblio}
\bibliographystyle{amsplainhyper}

\end{document}